\documentclass[12pt]{article}
\usepackage{setspace}
\usepackage{textcomp}
\usepackage{amsmath,amsthm,amssymb}
\usepackage[left=1.5in,top=1in,right=1in,nohead]{geometry}
\usepackage{latexsym}
\usepackage{mathrsfs}
\usepackage{rawfonts}
\usepackage[dvips]{graphicx}
\usepackage{stmaryrd}

\usepackage[sc,osf]{mathpazo}

\def\bct{\begin{center}}
\def\ect{\end{center}}
\def\beg{\begin}

\def\<{\langle}
\def\>{\rangle}

\def\mbb{\mathbb}
\def\mbbz{\mathbb Z}

\def\ol{\overline}

\def\tn{\textnormal}

\newtheorem{thm}{Theorem}[section]

\newtheorem{rmk}[thm]{Remark}

\title{Self-Dual metrics on non-simply connected 4-manifolds} \author{H\"ulya Arg\"uz   \and Mustafa Kalafat \and Y\i ld\i ray Ozan}

\begin{document}
\maketitle
\begin{abstract} We construct self-dual(SD) but not locally conformally flat(LCF) metrics 
on families of non-simply connected 4-manifolds with small signature. We construct various sequences with bounded or unbounded Betti numbers and Euler characteristic. These metrics have negative scalar curvature. As an application, this addresses the Remark 4.79 of \cite{besse}. 
  \end{abstract}



\section{Introduction}

\vspace{.05in}

Self-dual manifolds were introduced by Penrose in \cite{penrose}, 
and they were put on a firm mathematical foooting by Atiyah, Hitchin and Singer in \cite{ahs}. 
After the explicit constructions of Poon, LeBrun, Joyce, and the connected sum theorem of Donaldson and Friedman,  Taubes proved that every 4-manifold admits
an anti-self-dual metric after taking connected sum with
$k \ol{\mathbb{CP}}_2$, where $k$ is sufficiently large \cite{taubes}. Although this is a very useful theorem, the anti-self-dual metric here  is not explicit. The minimal number for $k$ is called the {\em Taubes Invariant}, which is unknown for 
most 
4-manifolds. The number $k$ in Taubes' work is very  large. 
There are many explicit constructions of self-dual metrics on simply connected manifolds, however there are very few examples for the non-simply connected case which are not locally conformally flat. The self-dual metrics
on the 
blow-ups of the 
connected sums of $S^3\times S^1$ \cite{clasdhblhopf,
kimsfasd}, and on the
blow-ups of 
$S^2\times \Sigma_g$ for $g\geq 1$ \cite{clasdblruled,klp
} are the only explicitly known examples. 
Here we give new examples of self-dual metrics on closed non-simply connected 4-manifolds, and show that many new topological types can be realized. 

\vspace{.05in}


The idea is first to construct LCF manifolds using the techniques introduced in \cite{handle}. See Figure \ref{sequenceofmetricssd1} for the construction of the manifolds. The manifolds in \cite{handle} do not satisfy the hypothesis for producing self-dual metrics with minimal number of blow-ups. So one needs to make modifications and produce a new family. One  needs to disconnect the boundary of the related 3-manifold. So that J.Kim's theorems are applicable and as a consequence these manifolds admit the hyperbolic ansatz self-dual metrics of LeBrun. We obtain the following results.



\begin{thm}\label{m5} 
The manifolds $M^5_{g,n} \sharp\, l\mathbb{CP}_2$ 
admit self-dual metric of negative scalar curvature for all $l\geq 2$. \end{thm}

\beg{cor}\label{bettigrowth} There are new infinite families of closed, non-simply connected,
self-dual 4-manifolds, with Betti number growth as follows:  
\beg{enumerate}
\item \label{SD-NCE1} $b_1\to\infty$, $b_2$ bounded, and $\chi\to-\infty$,
\item \label{SD-NCE2} $b_1\to\infty$, $b_2\to\infty$, and $\chi\to-\infty$,
\item $b_1\to\infty$, $b_2\to\infty$, and $\chi$ bounded, 
\item $b_1$ bounded, $b_2\to\infty$, and $\chi\to\infty$,
\item $b_1\to\infty$, $b_2\to\infty$, and $\chi\to\infty$.
\end{enumerate}
These manifolds have strictly negative scalar curvature. \end{cor}

We also 
prove the following theorem using Kim's result. Note that this cannot 
be obtained automatically from Taubes' theorem 
since it does not say anything about the scalar curvature. See \cite{handle} for the constructions of the manifolds in this theorem.

\begin{thm}\label{main} 
The manifolds $M^1_g \sharp\, l\mathbb{CP}_2$, $M^2_{g,n} \sharp\, l\mathbb{CP}_2$, $M^4_1 \sharp\, l\mathbb{CP}_2$ 
admit self-dual metrics of negative scalar curvature for all sufficiently large  $l$.   \end{thm}
Since their signature is nonzero, the above manifolds do not admit any locally conformally flat(LCF) metrics. In Remark 4.79 of \cite{besse}, 
examples of compact, self-dual (half-conformally flat) manifolds which are not conformally Einstein are asked for. The family of metrics in the parts \ref{SD-NCE1} and \ref{SD-NCE2} of Corollary \ref{bettigrowth} has negative Euler characteristic. This violates the Hitchin-Thorpe inequality; see \cite{hitchinthorpeinequality0} or the more recent \cite{hitchinthorpeinequality1,hitchinthorpeinequality2}. If there is an Einstein or conformally Einstein metric on these compact 4-manifolds then the Euler characteristic has to be non-negative by Hitchin-Thorpe. So  these metrics are instances in the remark which are explicit, non-LCF and of negative scalar curvature type.

\vspace{.1in}

{\bf Acknowledgements.} We would like to thank to Claude LeBrun for directing us to the field, and Jongsu Kim for encouragement.  The figure was constructed by using the IPE software of Otfried Cheong.

\newpage
\section{Background and Proofs}\label{sec1}

To enable the construction of the self-dual metrics, we first construct locally conformally flat(LCF) 
metrics using Braam's conformal compactification procedure given in \cite{braam}. One starts with a hyperbolic 3-manifold $N$ with boundary which is obtained from the hyperbolic space by taking the quotient with a cusp-free geometrically finite Kleinian group. Then one spins around the boundary to get a closed Riemannian 4-manifold X,
$$X:={N\times S^1 \over \{ b\times S^1 | \, b\in \partial N\} }.$$
Actually, this corresponds to crossing the 3-manifold with a circle and then contracting the circles that lie on the boundary. The above process coupled with the magnetic monopoles of LeBrun \cite{clexplicit} is called the {\em hyperbolic ansatz}, which 
yields self-dual metrics on the blow-ups. The following theorem of Jongsu Kim tells us the precise conditions needed to be able to carry out the hyperbolic ansatz.

\beg{thm}[\cite{kimsfasd}]  \label{kim}
If $N=\mathcal H^3/\Gamma$ is a noncompact hyperbolic 3-manifold obtained from a cusp-free 
geometrically finite Kleinian group $\Gamma$ and $H_2(\bar N;\mathbb R)/ H_2(\partial N;\mathbb R)$ is at most $1$-dimensional, then there exist self-dual metrics on $X\sharp\, l \mathbb{CP}_2$ 
for all sufficiently large natural number $l$ where $X$ is obtained from $N$ by spinning around its boundary. The sign of the conformal class of these metrics is identical to the sign of 
$1-d(\Lambda(\Gamma))$.

If moreover $H_2(\bar N;\mathbb R)=H_2(\partial N;\mathbb R)$ we can take $l$ to be bigger than or equal to the number of connected components of $\partial N$. \end{thm} 
\noindent 
To check the homological condition appearing in the hypothesis of the theorem, we need to use the following isomorphism \cite{braam} where the coefficients are integers: 
\begin{equation}
H_2(\Bar N)\oplus H_1(\Bar N,\partial \bar N)\tilde\longrightarrow H_2(X).\label{isomorphism}
\end{equation}
To establish this isomorphism one needs a little bit of work, details of which are not given in  \cite{braam} so we will give a proof here. Rather than shrinking the boundary circles, one can
attach 2-discs to fill out those circles; this is called the capping \cite{handle} operation. So one can work with the following decomposition of the 4-manifold:
$$X=\bar N\times S^1 \cup \partial\bar N \times D^2~~\textnormal{where}~~\bar N\times S^1 \cap \partial\bar N \times D^2=\partial\bar N \times S^1.$$
If we write down the Mayer-Vietoris exact sequence \cite{hatcher} for this pair we get the following piece:
$$\hspace{-8mm}\cdots\to H_2(\partial N\times S^1)\stackrel{A}{\to} H_2(N\times S^1)\oplus H_2(\partial N)
\stackrel{B}{\to} H_2(X)\stackrel{C}{\to} H_1(\partial N\times S^1)\stackrel{D}{\to} H_1(N\times S^1)\oplus H_1(\partial N)\to\cdots$$
To analyze this sequence we need the maps involved in the relative exact sequence:
$$\cdots\to H_1(\partial N)\stackrel{i_*^1}{\to} H_1(N)\stackrel{j_*^1}{\to} H_1(N,\partial N)\stackrel{\partial_*^1}{\to} H_0(\partial N) \stackrel{i_*^0}{\to} H_0(N)\to 
0\,\cdot$$
From this point on, under the assumption that $\partial N$ has one or two connected components, one gets an easier proof. Nonetheless we will prove the general assertion. 
We make up the following exact 
diagram.
$$\beg{array}{cc@{}ccc}
 && H_2(X) ~~\stackrel{C}{\longrightarrow}~~ \tn{Ker}\, i_*^0 && \\
 & \stackrel{B~~~}{\nearrow} && \searrow & \\   
H_2(\partial N)\oplus H_1(\partial N)\stackrel{A}{\longrightarrow} H_2(N)\oplus H_1(N)\,\oplus\, H_2(\partial N) &  &  
 && 0 \,\cdot \\    
 & \stackrel{\searrow}{\varphi~~~} && \nearrow & \\
 && H_2(N)\oplus\,\tn{Im}j_*^1   &&
  \end{array}$$
  
\noindent The domain and range of the map $A$ are decomposed into basic components 
according to the K\"unneth formula, and given by $A(a,b)=(i_*a,i_*b,-a)$. 
Using these conventions we define the new map $\varphi$ by $\varphi(x,y,z)=(x+i_*z\, , \, j_*^1 y)$. One can easily check the exactness by chasing through the relative exact sequence, where the map $j_*^1$ is also involved. So $\tn{Ker}\,\varphi=\tn{Im}A$ and  $\varphi$ is surjective.  
The upper part is taken from the Mayer-Vietoris sequence above, 
where the map $D$ is defined by $D=I\oplus -I$, and after decomposing into components  according to K\"unneth, one can make the identification $\tn{Ker}\, D=\tn{Ker}\, i_*^0$. Since  $\tn{Coker}A\approx \tn{Im}\,\varphi$, the short exact sequence with nucleus $H_2(X)$ in the Mayer-Vietoris sequence gives
$$0\longrightarrow H_2(N)\oplus\,\tn{Im}j_*^1 \stackrel{\bar B}{\longrightarrow} H_2(X) \stackrel{C}{\longrightarrow} \tn{Ker}\, i_*^0 \longrightarrow 0\,\cdot$$ 
On the other hand, since $\tn{Coker}\, i_*^1\approx \tn{Im}\, j_*^1$, pulling out the 
short exact sequence with nucleus $H_1(N,\partial N)$ from the relative sequence reads, 
$$0\longrightarrow \tn{Im}j_*^1 \stackrel{inc}{\longrightarrow} H_1(N,\partial N)\stackrel{\partial_*^1}\longrightarrow \tn{Ker}\, i_*^0 \longrightarrow 0\,\cdot$$ 
All the maps are natural. Now $H_0(\partial N)$ is always free, so these exact sequences split, and combining the two parts yields the isomorphism (\ref{isomorphism}). 

\vspace{.05in}

Next using the Lefschetz duality $H_1(N,\partial N)\approx H^2(N)$ and the universal coefficients theorem, after canceling the $\tn{Ext}$ terms,  the isomorphism becomes
$$\tn{Hom}(H_2(N),\mathbb Z)\oplus H_2(N) \approx \tn{Hom}(H_2(X),\mathbb Z).$$
Writing this according to the free and torsion pieces it takes the form
$$F_2N\oplus F_2N\oplus T_2N\approx F_2X.$$
So comparing the two sides we obtain\, $T_2N=0$ \, and \, $2F_2N\approx F_2X$.

\vspace{.05in}

\begin{figure}[!h] \bct
\includegraphics[width=\textwidth]{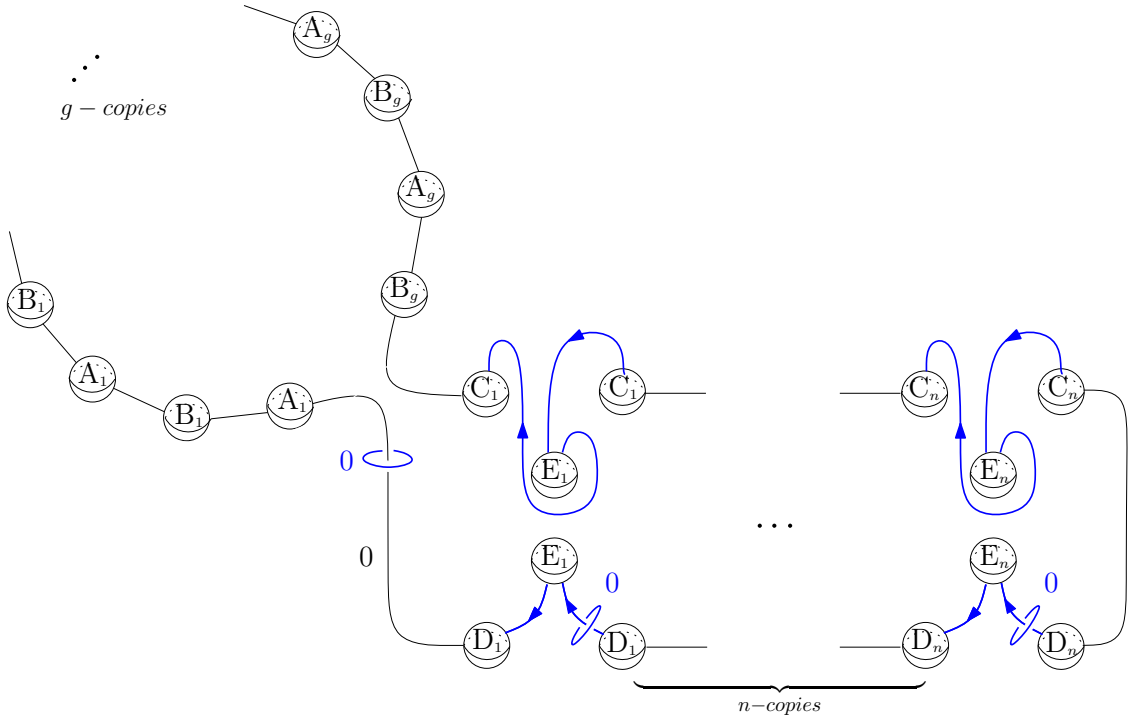}\\ \ect
  \caption{{\em {\small The LCF manifolds $M^5_{g,n}$.} } } \label{sequenceofmetricssd1}
\end{figure}

Now the paneled web 4-manifolds $M^1_g$, $M^2_{g,n}$
and $M_1^4$ of \cite{handle} have $b_2=2$, and by the above this gives $b_2(N)=1$,  which satisfies
the first hypothesis of Theorem \ref{kim}. 
However, they do not satisfy the second homological hypothesis. 
Next we construct a new set of manifolds to handle this situation. 
We modify the sequence $M^2_{g,n}$. To satisfy the second hypothesis,
 the boundary surfaces should generate the second homology of the $3$-manifold $N$. 
 The examples in \cite{handle} have all connected boundaries, and the boundary bounds
  the $3$-manifolds and is homologous to zero; hence the map 
  $$i_*: H_2(\partial N,\mbbz)\to H_2(N,\mbbz)$$
is zero. In order to get a nontrivial image, we should disconnect the boundary. 
So whenever we are identifying the boundary cylinders in the construction of the manifolds $M^2_{g,n}$ we identify in a parallel way, this time to get two distinct boundary components. In this way we obtain the sequence of manifolds $M^5_{g,n}$ which is shown concretely in Figure \ref{sequenceofmetricssd1}.
The following theorem ensures that the boundary surfaces generate $H_2(N,\mbbz)$.

\beg{thm} Let $N$ be a compact orientable $3$-manifold with two boundary components and $b_2=1$. 
Then each of the boundary surfaces is a generator of $H_2(N,\mbbz)$.
\end{thm}
\beg{proof} Take two copies of the manifold $N_1$ and $N_2$ such that the second copy is reversely oriented. Then identify them along their boundaries to form the compact $3$-manifold $M=N_1\cup_{\partial N} N_2$ without boundary, which we call {\em the double} of $N$. If we look at the Mayer-Vietoris exact sequence in this case, we see the following piece: 
$$\begin{array}{ccccc}
0\longrightarrow H_3(M)\longrightarrow & H_2(\partial N)& \stackrel{i^1_*\oplus -i^2_*}{\longrightarrow} & H_2(N_1)\oplus H_2(N_2) & \stackrel{j^1_*+j^2_*}{\longrightarrow} H_2(M)\longrightarrow \cdots\\
& \wr\wr           & &  \wr\wr\\   
& \mbbz\oplus\mbbz & & \mbbz\oplus\mbbz & \end{array}$$
\vspace{.05in}\\
For the second isomorphism here it can be checked through the universal coefficients theorem that there is no torsion. Assuming that the maps act through canonical isomorphisms, the images of the generators are 
$$i^1_*\oplus -i^2_* : (1,0)\mapsto (k,k), ~(0,1)\mapsto (l,l)$$
for some $k,l\in \mbbz$. Since $\partial N$ bounds both $N_1$ and $N_2$, the image 
of $(1,0)+(0,1)$ under $i^1_*$ and $i^2_*$ is zero implying $l=k$. Since by compactness $H_2(M)$ is torsion free, the isomorphism  following 
from exactness: $$ \mbbz\oplus\mbbz / \langle k(1,1)\rangle \approx H_2(M)$$
forces $k=\pm 1$, and hence it is a generator in both components. Furthermore, we have $H_2(M)=\mbbz$.\end{proof}

\noindent Since both of the homological conditions of \cite{kimsfasd} are satisfied, Theorem \ref{m5} follows.

\vspace{.05in}

Next we compute the topological invariants. The fundamental group
$$\llap{$\langle a$}_1
\cdots 
b_g, c_1\cdots c_n, d_1 \cdots d_n, e_1 \cdots e_n
~|~a_1^{-1}b_1^{-1}a_1b_1\cdots a_g^{-1}b_g^{-1}a_gb_gc_1\cdots c_n d_n^{-1}\cdots d_1^{-1},e_id_ie_i^{-1}c_i \rangle\qquad$$
gives relations $d_i=c_i^{-1}$ for $i=1\cdots n$\, and \, $(c_1\cdots c_n)^2=1$, after abelianization. Introducing a new variable $\bar c:=c_1\cdots c_n$ and dropping
the redundant  $c_n=c_{n-1}^{-1}\cdots c_1^{-1}\bar c$ ~yields
$$H_1(M^5_{g,n};\mbbz)=\langle a_1
\cdots 
b_g, c_1\cdots c_{n-1},\bar c, e_1 \cdots e_n
~|~ \bar c^2 \rangle_{+} ~\approx~ \mbbz^{2g+2n-1}\oplus \mbbz_2.$$
Counting the handles gives $\chi=
4-4g-4n$ \, and \, $b_2=0$, then \, $H_2(M^5_{g,n};\mbbz)=\mbbz_2.$
After taking the connected sum with $l$ copies of $\mbb{CP}_2$, we get the following invariants for the manifolds $M^5_{g,n} \sharp\, l\mathbb{CP}_2$ 
$$b_1=2g+2n-1, ~b_2=l, ~\chi=4-4g-4n+l.$$

Now, the sequences in Corollary \ref{bettigrowth} \, can be obtained by letting $g,n\longrightarrow\infty$, and taking $l=g+2$ or $l=4(g+n)$ in addition. Letting $l\longrightarrow\infty$, and taking $l=5(g+n)$ in addition. Alternatively taking $l=g+n$ one can again produce other sequences of the second type.

\begin{rmk} The sign of $1-d(\Lambda(\Gamma))$ is already computed to be negative 
 in \cite{handle} for the paneled web groups. This assures the sign of the scalar curvature in Theorem \ref{m5} as a consequence of Theorem \ref{kim} of Kim.
\end{rmk}

\vspace{.2in}

{\small \beg{flushleft} 
\textsc{
Universit\"at Hamburg, Bundesstrasse 55, 20146, Germany}\\
\textit{E-mail address:} \texttt{\textbf{huelya.arguez@\,math.uni-hamburg.de}} \end{flushleft}
}

{\small 
\beg{flushleft} \textsc{Tuncel\' \i  ~\" Un\' \i vers\' ites\' i, Turkia}\\
\textit{E-mail address:}  \texttt{\textbf{kalafg@\,gmail.com}} \end{flushleft}
}

{\small 
\beg{flushleft} \textsc{Orta Do\u gu Tekn\' \i k  \" Un\' \i vers\' ites\' i, 06800, Ankara, Turkia}\\
\textit{E-mail address:}  \texttt{\textbf{ozan@\,metu.edu.tr}} \end{flushleft}
}

\vspace{.05in}


\end{document}